\documentclass[12pt]{article}
\usepackage{latexsym, amssymb, amsmath, amscd, amsfonts, epsfig, graphicx, colordvi,verbatim,ifpdf}
\usepackage{amsfonts, amsmath, amssymb}
\usepackage{amssymb,amsfonts,amsmath,latexsym,epsfig,cite, psfrag,eepic,color}
\usepackage{amscd,graphics,supertabular}
\usepackage{latexsym, amssymb,  amsmath,amscd, amsfonts, epsfig, graphicx, colordvi,amsthm}

\usepackage{graphicx}
\usepackage{color}
\usepackage{ifpdf}
\usepackage{fancybox}
\usepackage{tikz}

\usepackage{float}
\usepackage{cases}

\usepackage{extarrows}

\newtheorem{thm}{Theorem}[section]

\newtheorem{defi}[thm]{Definition}

\def\pf{\noindent{\it Proof.} }
\def\qed{\nopagebreak\hfill{\rule{4pt}{7pt}}
\medbreak}

\makeatletter \@addtoreset{table}{section}
\numberwithin{equation}{section}

\def\qed{\nopagebreak\hfill{\rule{4pt}{7pt}}
\medbreak}

\newlength{\boxedparwidth}
  {\begin{center} \begin{tabular}{|@{\hspace{.315in}}c@{\hspace{.15in}}|}
                  \hline \\ \begin{minipage}[t]{\boxedparwidth}
                  \setlength{\parindent}{.25in}}%
  {\end{minipage} \\ \\ \hline \end{tabular} \end{center}}

\begin{document}
\begin{center}

 { \large\bf On the Enumeration and Congruences for $m$-ary Partitions}
\end{center}

\begin{center}
{Lisa Hui Sun}$^{1}$ and {Mingzhi Zhang}$^{2}$ \vskip 2mm

   Center for Combinatorics, LPMC\\
   Nankai University, Tianjin 300071, P. R. China\\[6pt]

   \vskip 2mm
    $^{1}$sunhui@nankai.edu.cn, $^{2}$mingzhi@mail.nankai.edu.cn
\end{center}

\vskip 6mm \noindent {\bf Abstract.} Let $m\ge 2$ be a fixed integer. Suppose that $n$ is a positive integer such that $m^j \leq n< m^{j+1}$ for some integer $j\ge 0$. Denote $b_{m}(n)$ the number of $m$-ary partitions of $n$, where each part of the partition is a power of $m$.  In this paper, we show that $b_m(n)$ can be represented as a $j$-fold summation by constructing a one-to-one correspondence between the $m$-ary partitions and a special class of integer sequences relying only on the base $m$ representation of $n$. It directly reduces to  Andrews, Fraenkel and Sellers' characterization of the values $b_{m}(mn)$ modulo $m$. Moreover, denote $c_{m}(n)$ the number of $m$-ary partitions of $n$ without gaps, wherein if $m^i$ is the largest part, then $m^k$ for each $0\leq k<i$ also appears as a part. We also obtain an enumeration formula for $c_m(n)$ which leads to an alternative representation  for the congruences of $c_m(mn)$ modulo $m$ due to Andrews, Fraenkel and Sellers.



\section{Introduction}
\allowdisplaybreaks

The arithmetic properties for partition functions have been extensively studied  since the discoveries of Ramanujan \cite{Ram19}. In this paper, we are mainly concerned with  the enumeration of $m$-ary partitions which leads to the congruence properties given by  Andrews, Fraenkel and Sellers \cite{And15, And16}.

Let $m\ge 2$ be a fixed  integer. An $m$-ary partition of a positive integer $n$ is a partition of $n$ such that each part is a power of $m$. The number of $m$-ary partitions of $n$ is denoted by $b_m(n)$.  For example, there are five $3$-ary partitions for $n=10$:
\begin{center}
9+1,\ 3+3+3+1,\ 3+3+1+1+1+1,\  3+1+1+1+1+1+1+1,\   1+1+1+1+1+1+1+1+1+1.
\end{center}
Thus, $b_3(10)=5$.  Denote an $m$-ary partition of $n$ by a sequence $\lambda=(a_\ell, a_{\ell-1}, \ldots, a_0)$ such that
\[
n=a_\ell m^\ell+a_{\ell-1} m^{\ell-1} +\cdots +a_0,
\]
where $a_\ell > 0$ and $a_i \geq 0$ for $0\leq i\leq \ell-1$. Denote the set of all the $m$-ary partitions of $n$ by $\mathcal{B}_m(n)$.
It is known that the generating function of $b_m(n)$ is given by
\begin{equation}\label{gf}
B_m(q)=\sum_{n=0}^\infty b_m(n)q^n =\prod_{k=0}^\infty \frac{1}{1-q^{m^k}}.
\end{equation}

For the case $m=2$, Churchhouse \cite{Chu69} conjectured the following congruences for the binary partition function $b_2(n)$:
\begin{align}\label{chu}
b_2(2^{2k+2}n)\equiv b_2(2^{2k}n) \pmod{2^{3k+2}},\nonumber\\[5pt]
b_2(2^{2k+1}n)\equiv b_2(2^{2k-1}n) \pmod{2^{3k}},
\end{align}
where $n,k\ge 1$. The conjecture was first proved by R{\o}dseth \cite{Rod70} and further studied by  Hirschhorn and Loxton \cite{HL75}. Later, it was extended to $m$-ary partitions by Andrews  \cite{And71}, Gupta \cite{Gup72}, and R{\o}dseth and Sellers \cite{RS01}.

Throughout this paper, without specification, we set $n$ to be a positive integer such that $m^j\leq n < m^{j+1}$ for some integer $j\ge 0$.
Recall that the base $m$ representation of $n$ is the unique expression of $n$ which can be written as follows
\begin{equation}\label{basem}
  n=\alpha_{j} m^{j}+\alpha_{j-1} m^{j-1}+\cdots+\alpha_{1} m+\alpha_{0},
\end{equation}
where $\alpha_j> 0$ and $0 \leq \alpha_{i}\leq m-1$ for $0\leq i \leq j-1$.  Denote the base $m$ representation of $n$ by
\begin{equation}\label{basem1}
r_{m}(n)=(\alpha_{j},\alpha_{j-1},\ldots,\alpha_{1},\alpha_{0}).
 \end{equation}
Based on the base $m$ representation of $n$, Andrews, Fraenkel and Sellers  \cite[Theorem 1]{And15}  provided the following modulo $m$ characterization of $b_m(mn)$:
\begin{equation}\label{modmc}
  b_{m}(mn)\equiv \prod\limits_{i=0}^{j}(\alpha_{i}+1) \pmod{m}.
\end{equation}

In this paper, by establishing a bijection between the set $\mathcal{B}_m(n)$ of $m$-ary partitions of $n$ and the set of  integer sequences given in the following theorem, we derive a $j$-fold summation formula for $b_m(n)$. It  will directly lead to Andrews, Fraenkel and Sellers' congruence \eqref{modmc}.

\begin{thm} \label{biject} There is a one-to-one correspondence between  the set $\mathcal{B}_m(n)$ of $m$-ary partitions of $n$ and the following set of integer sequences
\[
\mathcal{S}_m(n)=\{(\beta_{j}, \beta_{j-1}, \ldots, \beta_{1}) \,|\, 0\leq \beta_{j} \leq \alpha_{j}\ \mbox{and}\ 0\leq \beta_{t} \leq \alpha_{t}+m\beta_{t+1} \ \mbox{for}\  1 \leq t \leq j-1\}.
\]
 \end{thm}

Based on the above bijection, we provide a combinatorial approach to derive the following $j$-fold summation formula for $b_m(n)$.

\begin{thm} \label{nmaryp} Let $r_{m}(n)$=$(\alpha_{j},\alpha_{j-1},\ldots,\alpha_{1},\alpha_{0})$ be the base $m$ representation of $n$. We have
\begin{equation}\label{enumermn}
  b_{m}(n)=\sum_{k_j=0}^{\alpha_{j}}\sum_{k_{j-1}=0}^{\alpha_{j-1}+mk_j}
  \ldots\sum_{k_1=0}^{\alpha_{1}+mk_2} 1.
\end{equation}
Obviously, $b_m(n)=1$ when $j=0$.
\end{thm}

Notice that if $r_{m}(n)$=$(\alpha_{j},\alpha_{j-1},\cdots,\alpha_{1},\alpha_{0})$, then  $r_{m}(mn)$=$(\alpha_{j},\alpha_{j-1},\cdots,\alpha_{1},\alpha_{0}, 0)$. Thus the above theorem leads to that
\[
 b_{m}(mn)=\sum_{k_j=0}^{\alpha_{j}}\sum_{k_{j-1}=0}^{\alpha_{j-1}+mk_j}
  \ldots\sum_{k_1=0}^{\alpha_{1}+mk_2}\sum_{k_0=0}^{\alpha_{0}+mk_1} 1.
\]
By taking modulo $m$ on both sides of the above equation, it directly reduces to Andrews, Fraenkel and Sellers' congruence \eqref{modmc}.

We also consider the cases for the $m$-ary partitions  without gaps,  wherein if $m^i$ is the largest part, then $m^k$ for each $0\leq k<i$ also appears as a part. The related works on such restricted  $m$-ary partitions can be found in \cite{And17, And16, Hou15}. Moreover,  in \cite{AS07, FHRT16, HS04, S05}, a general class of non-squashing partitions was introduced and studied, which contains  $m$-ary partitions as a special case.

Let $c_{m}(n)$ denote the number of $m$-ary partitions without gaps of $n$. Based on the bijection given in Theorem \ref{biject}, we also obtain the following enumeration formula for $c_m(n)$.

\begin{thm} \label{nmarypnogap}
Let $r_{m}(n)=(\alpha_{j},\alpha_{j-1},\ldots,\alpha_{1},\alpha_{0})$ be the base $m$ representation of $n$. We have
\begin{equation}\label{enumercmn}
  c_{m}(n)=1+\sum_{r=1}^{j}\sum_{k_{r}=\chi_{r}}^{\lfloor\frac{n}{m^{r}}\rfloor-1}
  \ldots\sum_{k_1=\chi_{1}}^{\alpha_{1}-1+mk_2} 1,
\end{equation}
where for $1\leq i \leq r$,
\begin{equation}\label{chi}
\chi_{i}=\left\{
           \begin{array}{ll}
             0, & \hbox{if $\alpha_{i-1} > 0$,} \\[5pt]
             1, & \hbox{if $\alpha_{i-1} = 0$.}
           \end{array}
         \right.
\end{equation}
\end{thm}

Applying formula \eqref{enumercmn}, we obtain the following congruence  property of $c_m(mn)$, which reveals the results given by Andrews, Fraenkel and Sellers \cite[Theorem 2.1]{And16}.

\begin{thm} Let $\chi_i$ be defined by \eqref{chi} for $1\leq i\leq r$, then we have
\begin{align}\label{concmn}
  c_{m}(mn)
  &\equiv \alpha_0+(\alpha_0-1)\sum_{i=1}^{j}(\alpha_1-\chi_1)(\alpha_2-\chi_2)\cdots(\alpha_i-\chi_i) \pmod m.
\end{align}
\end{thm}

\section{The  enumeration formula for $b_m(n)$}

 In this section, we provide a bijection between the set $\mathcal{B}_m(n)$ of $m$-ary partitions of $n$ and the set
\[
\mathcal{S}_m(n)=\{(\beta_{j}, \beta_{j-1}, \ldots, \beta_{1}) \,|\, 0\leq \beta_{j} \leq \alpha_{j}\ \mbox{and}\ 0\leq \beta_{t} \leq \alpha_{t}+m\beta_{t+1} \ \mbox{for}\  1 \leq t \leq j-1\},
\]
 which relies only on the base $m$ representation of $n$. It will lead to the enumeration formula \eqref{enumermn} for the  $m$-ary partitions.

To this end, we first define the following subtraction between the base $m$ representation of $n$ and an ordinary $m$-ary partition of $n$.

\begin{defi}\label{def-subs} Let $r_m(n)=(\alpha_j, \alpha_{j-1}, \ldots, \alpha_0)$ be the base $m$ representation of $n$ and $\lambda=(\lambda_{\ell},\lambda_{\ell-1},\ldots, \lambda_{0})$ be an $m$-ary partition of $n$. Then subtracting $\lambda$ from $r_m(n)$ is given as follows
\begin{equation}\label{subt}
r_m(n) -\lambda = (
 \beta_j,
\beta_{j-1},\ldots,\beta_1),
\end{equation}
where for $1\leq i \leq j$,
\[
\beta_{i}=\sum_{k=i}^{j} m^{k-i}(\alpha_{k}-\lambda_{k})
\]
provided that $\lambda_k=0$ for  $\ell <k\leq j$.

\end{defi}

We further show that the  subtraction \eqref{subt} defined above gives a bijection between $\mathcal{B}_m(n)$ and $\mathcal{S}_m(n)$.

\begin{thm}\label{thmbij}
Let $r_m(n)=(\alpha_j, \alpha_{j-1}, \ldots, \alpha_0)$ be the base $m$ representation of $n$ and $\lambda=(\lambda_{\ell},\lambda_{\ell-1},\ldots, \lambda_{0})$ be an arbitrary $m$-ary partition of $n$.
Define a map $\varphi$ from  $\mathcal{B}_m(n)$ to $\mathcal{S}_m(n)$  by $\varphi(\lambda)=r_m(n) - \lambda$.
 Then $\varphi$ is a bijection between   $\mathcal{B}_m(n)$  and $\mathcal{S}_m(n)$.
\end{thm}

\pf  Denote $\beta=\varphi(\lambda)=(\beta_j,\beta_{j-1}, \ldots, \beta_1)$. First, we proceed to show that $\beta\in \mathcal{S}_m(n)$ and thereby $\varphi$ is well defined.
Following Definition \ref{def-subs}, it is easy to see that
\begin{align}
&\beta_j=\alpha_j-\lambda_j,\label{betaj}\\[5pt]
&\beta_{t}=\alpha_{t}-\lambda_{t}+m \beta_{t+1},\label{betat}
\end{align}
where $1\leq t \leq j-1$ and $\lambda_k=0$ for  $\ell <k\leq j$.  Since $\lambda_k\ge 0$ for  $0 \leq k\leq j$,  we see that $ \beta_j \leq \alpha_j$ and $\beta_{t}\leq \alpha_{t}+m \beta_{t+1}$ for $1\leq t \leq j-1$.

It is obvious that $\lambda_j\leq \alpha_j$, so that $\beta_j\ge 0$. From the fact that
\[
\lambda_{j} m^{j}+\lambda_{j-1} m^{j-1}+\cdots+\lambda_{0}=\alpha_{j} m^{j}+\alpha_{j-1} m^{j-1}+\cdots+\alpha_{0},
\]
we are led to that for $1\leq t \leq j-1$,
\[
\lambda_{j} m^{j}+\lambda_{j-1} m^{j-1}+\cdots+\lambda_{t}m^{t} \leq \alpha_{j} m^{j}+\alpha_{j-1} m^{j-1}+\cdots+\alpha_{0}.
\]
Hence we obtain that
\[
\Big((\lambda_{t}-\alpha_{t})+\sum_{k=1}^{j-t} (\lambda_{t+k}-\alpha_{t+k})m^k\Big) m^{t} \leq
\alpha_{t-1} m^{t-1}+\alpha_{t-2} m^{t-2}+\cdots+\alpha_{0}.
\]
Since $(\alpha_j, \alpha_{j-1},\ldots, \alpha_0)$ is the base $m$ representation of $n$, it is obvious that
\[
\alpha_{t-1} m^{t-1}+\alpha_{t-2} m^{t-2}+\cdots+\alpha_{0}< m^{t},
\]
which implies that
\[
\Big((\lambda_{t}-\alpha_{t})+\sum_{k=1}^{j-t} (\lambda_{t+k}-\alpha_{t+k})m^k\Big) m^{t} < m^{t}.
\]
Note that $\lambda_k$ and $\alpha_k$ are all integers for $0\leq k\leq j$, it follows that
\[
\lambda_{t}-\alpha_{t}+\sum_{k=1}^{j-t} (\lambda_{t+k}-\alpha_{t+k})m^k\leq 0
\]
and thereby
\begin{align*}
&\lambda_{t} \leq \alpha_{t}+\sum_{k=1}^{j-t} (\alpha_{t+k}-\lambda_{t+k})m^k = \alpha_{t}+m \beta_{t+1}.
\end{align*}
By \eqref{betat}, it directly leads to that  $\beta_{t}\geq 0$ for $1\leq t\leq j-1$. Thus $\beta \in \mathcal{S}_m(n)$ and $\varphi$ is well defined.

To prove that $\varphi$ is a bijection, it is sufficient to show that there exists the inverse map of $\varphi$. For a given  $\beta=(\beta_{j}, \beta_{j-1}, \ldots, \beta_{1}) \in \mathcal{S}_m(n)$, let  $\varphi^{-1}(\beta)$ be given by computing
\[
\lambda'=(\alpha_{j}-\beta_{j},\,\alpha_{j-1}-\beta_{j-1}+m\beta_{j},\,\ldots,\, \alpha_{1}-\beta_{1}+m\beta_{2},\,\alpha_{0}+m\beta_{1})
 \]
 and then  deleting the preceding zeros.
From the definition of $\mathcal{S}_m(n)$, we see that each element of $\lambda'$ is nonnegative. Furthermore, it is easy to see that
\begin{align*}
(\alpha_{j}&-\beta_{j})m^j+(\alpha_{j-1}-\beta_{j-1}+m\beta_{j}) m^{j-1}+\cdots+\alpha_{0}+m\beta_{1}=n,
\end{align*}
which implies that $\varphi^{-1}(\beta)$  is an $m$-ary partition of $n$ and thereby $\varphi^{-1}(\beta) \in \mathcal{B}_m(n)$.
It completes the proof of the bijection. \qed

For example, let  $m=4$ and $n=36$, then the base $4$ representation of $36$ is $r_4(36)=(2,1,0)$.  The correspondence between
all the $4$-ary partitions of $36$ and the integer sequences belonging to $\mathcal{S}_4(36)$ can be seen in Table~2.1.
\begin{table}[h]
\caption{The correspondence between $\lambda \in \mathcal{B}_4(36)$ and $\beta \in \mathcal{S}_4(36)$}
\begin{center}
{\small
  \begin{tabular}{cccccccccc}
    \hline
    $\lambda$ &$(2,1,0)$&$(2,0,4)$&$(1,5,0)$&$(1,4,4)$&$(1,3,8)$&$(1,2,12)$&$(1,1,16)$&$(1,0,20)$&$(0,9,0)$
    \\
    $\beta$ &$(0,0)$&$(0,1)$&$(1,0)$&$(1,1)$&$(1,2)$&$(1,3)$&$(1,4)$&$(1,5)$&$(2,0)$\\  \hline

    $\lambda$ &$(0,8,4)$&$(0,7,8)$&$(0,6,12)$&$(0,5,16)$&$(0,4,20)$&$(0,3,24)$&$(0,2,28)$&$(0,1,32)$&$(0,0,36)$\\

    $\beta$ &$(2,1)$&$(2,2)$&$(2,3)$&$(2,4)$&$(2,5)$&$(2,6)$&$(2,7)$&$(2,8)$&$(2,9)$\\   \hline
  \end{tabular}\label{table1}
}
  \end{center}
\end{table}

The above theorem directly leads to that $b_m(n)=|\mathcal{B}_m(n)|=|\mathcal{S}_m(n)|$. By studying the recursive properties of the sequences in $\mathcal{S}_m(n)$, we obtain the  $j$-fold summation  formula \eqref{enumermn} of $b_m(n)$. Now we give the detailed proof of Theorem \ref{nmaryp}.

{\noindent \it Proof of Theorem \ref{nmaryp}.}
 Denote the summation on the right hand side of \eqref{enumermn} by
 \[
 f(\alpha_{j},\alpha_{j-1},\ldots,\alpha_{1},\alpha_{0})=\sum_{k_j=0}^{\alpha_{j}}\sum_{k_{j-1}=0}^{\alpha_{j-1}+mk_j}
  \ldots\sum_{k_2=0}^{\alpha_{2}+mk_3}\sum_{k_1=0}^{\alpha_{1}+mk_2} 1.
 \]
 We prove the theorem by induction. When $j=0$,  $r_m(n)=(\alpha_0)$ with $0 \leq \alpha_0 \leq m-1$. It is obvious that $b_m(n)=f(\alpha_0)=1$.

Suppose that \eqref{enumermn} holds for $j=i-1$. When $j=i$, we have
$r_m(n)=(\alpha_{i},\alpha_{i-1},\ldots,\alpha_{1},\alpha_{0})$.
   By Theorem \ref{thmbij}, it implies that $b_m(n)=|\mathcal{S}_m(n)|$ where
\[
\mathcal{S}_m(n)=\{(\beta_{i}, \beta_{i-1},\ldots, \beta_1) \,|\, 0\leq \beta_{i} \leq \alpha_{i}, 0\leq \beta_{t} \leq \alpha_{t}+m\beta_{t+1}, 1 \leq t \leq i-1 \}.
\]
For a fixed $\beta_{i}$ with $0\leq \beta_{i} \leq \alpha_{i}$, let us consider the subset of $\mathcal{S}_m(n)$  with $\beta_{i}$ being the first entry. By deleting $\beta_{i}$ in these sequences, it is easy to see that this subset is bijective with the following set
\[
\mathcal{S}_{\beta_{i}}=\big\{(\beta_{i-1}, \beta_{i-2},\ldots, \beta_1) \,|\,  0\leq \beta_{t} \leq \alpha_{t}+m\beta_{t+1}, 1 \leq t \leq i-1 \big\},
\]
and therefore
\[
\mathcal{S}_m(n)=\bigcup_{\beta_i=0}^{\alpha_i} \mathcal{S}_{\beta_i}.
\]
By induction, we obtain that the cardinality of the set $\mathcal{S}_{\beta_{i}}$ is
\begin{align*}
 |\mathcal{S}_{\beta_{i}}|=f(\alpha_{i-1}+m \beta_{i}, \alpha_{i-1}, \ldots, \alpha_0)=\sum_{k_{i-1}=0}^{\alpha_{i-1}+m \beta_{i}}\ldots\sum_{k_2=0}^{\alpha_{2}+mk_3}\sum_{k_1=0}^{\alpha_{1}+mk_2} 1.
\end{align*}
Then by summing the above equation for $\beta_{i}$ from $0$ to $\alpha_{i}$, we obtain
\begin{align*}
 |\mathcal{S}_m(n)|=f(\alpha_{i}, \alpha_{i-1}, \ldots, \alpha_0)=\sum_{k_{i}=0}^{\alpha_{i}}\sum_{k_{i-1}=0}^{\alpha_{i-1}+mk_i}\ldots\sum_{k_2=0}^{\alpha_{2}+mk_3}\sum_{k_1=0}^{\alpha_{1}+mk_2} 1,
\end{align*}
which completes the proof. \qed

Note that the  $j$-fold summation formula \eqref{enumermn} also can be derived from the generating function \eqref{gf} of $b_m(n)$. Moreover, by setting $M=\{m,m,\ldots\}$ in the summation given by Folsom, Homma,  Ryu and Tong \cite[Theorem 1.5]{FHRT16}, it reduces to another $j$-fold summation expression for $b_m(n)$.

\section{The $m$-ary partitions without gaps}

In this section, based on the  bijection given in Theorem \ref{thmbij}, we derive an enumeration formula for the $m$-ary partitions without gaps. Denote $c_m(n)$ the number of this restricted $m$-ary partitions of $n$. We also obtain an alternative expression for the congruence properties of $c_m(mn)$ given by  Andrews, Fraenkel and Sellers \cite[Theorem 2.1]{And16}.

Recall that by using the base $m$ representation of $n$ in the following form
\[
n=\sum_{i=\ell}^{\infty}\alpha_{i}m^{i}
\]
where $1 \leq \alpha_{\ell} < m$ and
$0 \leq \alpha_{i} <m$ for $i > \ell$, Andrews, Fraenkel and Sellers obtained the following result.

\begin{thm}[Andrews, Fraenkel and Sellers {\cite[Theorem 2.1]{And16}}] \label{ThmAnd16}

\begin{itemize}
  \item[{\rm (1)}] If $\ell$ is even, then
\begin{equation}\label{And-1}
c_{m}(mn) \equiv \alpha_{\ell}+(\alpha_{\ell}-1)\sum_{i=\ell+1}^{\infty}\alpha_{\ell+1}\cdots \alpha_{i} \pmod m.
\end{equation}
  \item[{\rm (2)}] If $\ell$ is odd, then
\begin{equation}\label{And-2}
c_{m}(mn) \equiv 1-\alpha_{\ell}-(\alpha_{\ell}-1)\sum_{i=\ell+1}^{\infty}\alpha_{\ell+1}\cdots \alpha_{i} \pmod m.
\end{equation}
\end{itemize}
\end{thm}

Denote the floor function of a real number $a$ by $\lfloor a \rfloor$, which is the largest integer less than or equal to $a$.
To derive our expression of the congruences \eqref{And-1} and \eqref{And-2}, first let us show how to derive the enumeration formula \eqref{enumercmn} for $c_m(n)$ as given in Theorem  \ref{nmarypnogap}.

{\noindent \it Proof of Theorem \ref{nmarypnogap}.}  Denote the set of all the $m$-ary partitions without gaps of $n$ by $\mathcal{G}_m(n)$.
We claim that for any $\lambda \in \mathcal{G}_m(n)$, it can be written as
\[
\lambda=\Big(\left\lfloor\frac{n}{m^{r}}\right\rfloor-\beta_{r},\, \alpha_{r-1}-\beta_{r-1}+m\beta_{r},\,
\ldots,\, \alpha_{0}+m\beta_{1}\Big),
\]
where $0 \leq r \leq j$ and $\beta_i$ are integers such that
\begin{equation}\label{greaterzero}
\left\lfloor\frac{n}{m^{r}}\right\rfloor-\beta_{r}>0,\
\alpha_{r-1}-\beta_{r-1}+m\beta_{r}>0,\
\ldots,\
\alpha_{0}+m\beta_{1}>0.
\end{equation}
Specially, when $r=0$, there is a unique $m$-ary partition without gaps, say, $\lambda=(n)$ which consists of $n$ ones.
For $r>1$, as an example, we consider $m=4$  and $n=73$, then $r_4(73)=(1, 0, 2, 1)$. When $r=2$,  we can obtain that $(3,6,1)$ is a $4$-ary partition without gaps which can be represented as $\big(\left\lfloor\frac{73}{4^2}\right\rfloor-1, 2-0+4\times 1, 1+4\times 0\big)$.

Recall that
\[
n=\alpha_{j} m^{j}+\alpha_{j-1} m^{j-1}+\cdots+\alpha_{0}.
\]
For $1\leq r \leq j$, it follows that
\begin{align*}
&\Big(\left\lfloor\frac{n}{m^{r}}\right\rfloor-\beta_{r}\Big)  m^{r}+(\alpha_{r-1}-\beta_{r-1}+m\beta_{r}) m^{r-1}+\cdots+(\alpha_{1}-\beta_{1}+m\beta_{2}) m + \alpha_{0}+m\beta_{1}\\[5pt]
&\quad=\left\lfloor\frac{n}{m^{r}}\right\rfloor  m^{r}+\alpha_{r-1} m^{r-1}+\cdots+\alpha_{1}  m + \alpha_{0}\\[5pt]
&\quad=(\alpha_{j} m^{j-r}+\cdots+\alpha_{r}) m^{r}+\alpha_{r-1} m^{r-1}+\cdots+\alpha_{1}  m + \alpha_{0}\\[5pt]
&\quad=\alpha_{j} m^{j}+\cdots+\alpha_{r} m^{r}+\alpha_{r-1} m^{r-1}+\cdots+\alpha_{1}  m + \alpha_{0}\\
&\quad=n,
\end{align*}
which certifies that $\lambda \in \mathcal{G}_m(n)$. For a fixed integer $r$ such that $0\leq r \leq j$,
by the bijection given in Theorem \ref{thmbij}, we get
\[
\varphi(\lambda)=\Big(\alpha_j,\, \alpha_{j-1}+m\alpha_j,\,\ldots, \,\sum_{k=r+1}^j m^{k-r-1}\alpha_k, \,\beta_r,\, \beta_{r-1},\, \ldots, \,\beta_1\Big).
\]
Note that for any $\lambda \in \mathcal{G}_m(n)$ with the given $r$, the first $j-r$ elements in $\varphi(\lambda)$ are the same, which only depend on $\alpha_{r+1}, \ldots, \alpha_j$. Then by deleting these terms, we find that the set of $m$-ary partitions without gaps $\mathcal{G}_m(n)$ is in one-to-one correspondence with the following set of integer sequences
\[
\mathcal{R}_m(n)=\bigcup_{r=0}^{j}\left\{(\beta_{r},  \ldots, \beta_{1})\,\big |\,  \left\lfloor\frac{n}{m^{r}}\right\rfloor-\beta_{r}>0,\,
\alpha_{r-1}-\beta_{r-1}+m\beta_{r}>0,\,
\ldots,\,
\alpha_{0}+m\beta_{1}>0\right\}.
\]
It indicates that $c_m(n)=|\mathcal{R}_m(n)|$.

From the conditions \eqref{greaterzero}, we see that if $\alpha_{0}=0$, then $\beta_{1}>0$, which means $\beta_{1}$ starts from 1. If $\alpha_{0}>0$, then $\beta_{1}\geq 0$, which means $\beta_{1}$ starts from 0. For both cases, we denote $\beta_{1}$ starting  from $\chi_{1}$, which is defined by \eqref{chi}. By $\alpha_{1}-\beta_{1}+m\beta_{2}>0$ we have $\beta_{1}<\alpha_{1}+m\beta_{2}$. Thereby we see that $\beta_{1}$  ranges from $\chi_{1}$ to $\alpha_{1}+m\beta_{2}-1$. By similar arguments applying to $\beta_{i}$ for $2 \leq i \leq r$, we have
\[
c_m(n)=|\mathcal{R}_m(n)|=1+\sum_{r=1}^{j}\sum_{k_{r}=\chi_{r}}^{\lfloor\frac{n}{m^{r}}\rfloor-1}
  \ldots\sum_{k_1=\chi_{1}}^{\alpha_{1}-1+mk_2} 1,
\]
where for $1\leq i \leq j$,
\begin{equation*}
\chi_{i}=\left\{
           \begin{array}{ll}
             0, & \hbox{if $\alpha_{i-1} > 0$,} \\[5pt]
             1, & \hbox{if $\alpha_{i-1} = 0$.}
           \end{array}
         \right.
\end{equation*}
This completes the proof. \qed

Noting that $r_m(mn)=(\alpha_j, \alpha_{j-1}, \ldots, \alpha_0, 0)$, then by applying the above result  we have
\[
c_m(mn)=1+\sum_{r=1}^{j+1}\sum_{k_{r}=\chi_{r-1}}^{\lfloor\frac{mn}{m^{r}}\rfloor-1}
  \ldots\sum_{k_2=\chi_1}^{\alpha_{1}-1+mk_3}\sum_{k_1=1}^{\alpha_{0}-1+mk_2} 1.
\]
By taking modulo $m$ on both sides of the above identity, we directly obtain the congruence property \eqref{concmn}, namely,
\begin{align}\label{concmn1}
  c_{m}(mn)
  &\equiv 1+\sum_{r=1}^{j+1}(\alpha_0-1)(\alpha_1-\chi_1)(\alpha_2-\chi_2)\cdots(\alpha_{r-1}-\chi_{r-1}) \pmod m \nonumber \\
&\equiv 1+(\alpha_0-1)+(\alpha_0-1)\sum_{r=2}^{j+1}(\alpha_1-\chi_1)(\alpha_2-\chi_2)\cdots(\alpha_{r-1}-\chi_{r-1}) \pmod m \nonumber \\
&\equiv \alpha_0+(\alpha_0-1)\sum_{i=1}^{j}(\alpha_1-\chi_1)(\alpha_2-\chi_2)\cdots(\alpha_{i}-\chi_{i}) \pmod m.
\end{align}

As an example, let $m=5$ and $n=485=3\cdot 5^3+4\cdot 5^2+2\cdot 5$. Then $j=3$ and the base $5$ representation of $485$ is $r_5(485)=(\alpha_3,\alpha_2,\alpha_1,\alpha_0)=(3,4,2,0)$.
Therefore
\[
\chi_3=\chi_2=0,\  \chi_1=1,
\]
and
\begin{align*}
c_5(5\cdot 485)&\equiv  -\big((\alpha_1-1)+(\alpha_1-1)\alpha_2+(\alpha_1-1)\alpha_2\alpha_3\big) \pmod{5}\\[5pt]
& =-(1+1\cdot 4+ 1\cdot 4\cdot 3)\\[5pt]
&=-17\equiv 3 \pmod{5}.
\end{align*}
In fact, we have $c_5(5\cdot 485)=230358 \equiv 3 \pmod{5}$, which coincides with the above result.

To conclude this paper, we remark that the congruence \eqref{concmn1} for $c_m(mn)$ is equivalent to Theorem \ref{ThmAnd16}   due to  Andrews, Fraenkel and Sellers \cite{And16}.

{\noindent \it Proof of Theorem \ref{ThmAnd16}.} Let $r_{m}(n)$=$(\alpha_{j},\alpha_{j-1},\ldots,\alpha_{1},\alpha_{0})$ be the base $m$ representation of $n$.

Following Lemma 2.9 of \cite{And16}, we see that $c_{m}(m^3n) \equiv c_{m}(mn) \pmod m$ for all $n \geq 0$. Thereby to prove Theorem \ref{ThmAnd16}, it is sufficient to show the cases that $\ell=0$ and $\ell=1$, which correspond to $\alpha_{0}>0$ and $\alpha_{0}=0$ ($\alpha_1 > 0$), respectively.

If $\alpha_{0} > 0$, then $\chi_1=0$. It leads to that
\begin{align}\label{concmn2}
c_{m}(mn)&\equiv \alpha_{0}+(\alpha_{0}-1)\sum_{i=1}^{j}\alpha_1(\alpha_2-\chi_2)\cdots(\alpha_i-\chi_i)\pmod{m}.
\end{align}
We further consider the values of $\alpha_1, \alpha_2, \cdots, \alpha_j$. If $\alpha_i > 0$ for $i \geq 1$, then $\chi_i=0$ for
$2 \leq i \leq j$. Thus \eqref{concmn2} turns to be
\begin{align*}
c_{m}(mn)&\equiv \alpha_{0}+(\alpha_{0}-1)\sum_{i=1}^{j}\alpha_1\cdots\alpha_i \equiv \alpha_{0}+(\alpha_{0}-1)\sum_{i=1}^{\infty}\alpha_1\cdots\alpha_i\pmod{m},
\end{align*}
where $\alpha_i=0$ for $i > j$. Otherwise, suppose that $\alpha_k \ (1 \leq k \leq j)$ is the first zero in the sequence $\alpha_1, \alpha_2, \cdots,$ then $\chi_i=0$ for $1 \leq i \leq k$. Noting that $\alpha_k=0$, we obtain
\begin{align*}
  \sum_{i=1}^{j}\alpha_1(\alpha_2-\chi_2)\cdots(\alpha_i-\chi_i) &=\sum_{i=1}^{k-1}\alpha_1\alpha_2\cdots\alpha_i =\sum_{i=1}^{\infty}\alpha_1\alpha_2\cdots\alpha_i,
\end{align*}
where $\alpha_i=0$ for $i > j$. Therefore, \eqref{concmn2} leads to that
\begin{align*}
c_{m}(mn) \equiv \alpha_{0}+(\alpha_{0}-1)\sum_{i=1}^{\infty}\alpha_1\cdots\alpha_i \pmod{m},
\end{align*}
and both cases coincide with \eqref{And-1} with $\ell=0$.

If $\alpha_{0} = 0$ and $\alpha_1 > 0$, following the same procedure, we obtain that
\begin{align*}
c_{m}(mn)&\equiv (-1)\sum_{i=1}^{j}(\alpha_1-1)\alpha_2(\alpha_3-\chi_3)\cdots(\alpha_i-\chi_i)\pmod{m}\\
&\equiv 1-\alpha_{1}-(\alpha_{1}-1)\sum_{i=2}^{j}\alpha_{2}(\alpha_3-\chi_3)\cdots(\alpha_i-\chi_i) \pmod{m}\\
&\equiv 1-\alpha_{1}-(\alpha_{1}-1)\sum_{i=2}^{\infty}\alpha_{2}\cdots\alpha_{i} \pmod{m},
\end{align*}
which coincides with \eqref{And-2} with $\ell=1$. This completes the proof. \qed

\vskip 15pt \noindent {\small {\bf Acknowledgments.}

This work was supported by the National Science Foundation of China, and the Natural Science Foundation of Tianjin, China.

\end{document}